\documentstyle[12pt]{article}
\def\Bbb#1{{\bf #1}}
\newtheorem{proposition}{Proposition}[section]
\newtheorem{lemma}[proposition]{Lemma}

\newtheorem{theorem}[proposition]{Theorem}
\begin{document}
\title{Adiabatic Limits and Spectral Geometry
of Foliations}
\author{Yuri A. Kordyukov\\Department of Mathematics,\\
Ufa State Aviation Technical University,\\
12 K.Marx str.,\\ Ufa, 450025, Russia\thanks{Current address (till August 31,
1995): Max-Planck-Institut
fuer Mathematik, Gottfried-Claren-Strasse 26,
53225 Bonn, Germany.
E-mail: kordy@mpim-bonn.mpg.de}}
\date{}
\maketitle
\begin{abstract}
We study spectral asymptotics for the Laplace operator
 on differential forms on a Riemannian foliated manifold equipped with a
bundle-like metric in the case when the
metric is blown up in directions normal to the leaves of
the foliation. The asymptotical
formula for the eigenvalue distribution function is obtained.
The relationships with the spectral theory of leafwise
Laplacian and with the noncommutative
spectral geometry of foliations are discussed.
\end{abstract}
{\bf Introduction}
\medskip
\par
\noindent Let $(M,{\cal F})$ be a closed foliated manifold,
$\dim M = n$, $\dim {\cal F} = p$, $p+q=n$,
equipped with a Riemannian metric $g_M$. We assume that
the foliation ${\cal F}$ is Riemannian, and the metric
$g_M$ is bundle-like.
Let $F=T{\cal F}$ be an integrable distribution
of $p$-planes in $TM$, and $H=F^{\bot}$ be the orthogonal
complement to $F$. So we have a decomposition of
$TM$ into a direct sum:
\begin{equation}
\label{decomp}
TM=F\bigoplus H.
\end{equation}
The decomposition~(\ref{decomp}) induces the decomposition
of the metric
\begin{equation}
g_{M}=g_{F}+g_{H}.
\end{equation}
Define a one-parameter family $g_{h}$ of
metrics on $M$ by the formula
\begin{equation}
g_{h}=g_{F} + {h}^{-2}g_{H},
0 < h \leq 1.
\end{equation}
For any $h>0$, we have the Laplace operator on differential forms defined by
the metric $g_h$:
\begin{equation}
\Delta_{h}=d^{*}_{g_h}d+dd^{*}_{g_h},
\end{equation}
where $d$ is the de Rham differential:
\begin{equation}
d:C^{\infty}(M,\Lambda^k T^{*}M)\rightarrow
C^{\infty}(M,\Lambda^{k+1}T^{*}M),
\end{equation}
$d^{*}_{g_h}$ is the adjoint with respect to
the metric on $C^{\infty}(M,\Lambda T^{*}M)$ induced by $g_{h}$.
The operator $\Delta_{h}$ is a self-adjoint, elliptic differential operator
with
the positive definite, scalar principal symbol
in the Hilbert space $L^2(M,\Lambda T^{*}M,g_h)$.
By the standard perturbation theory, there are (countably many)
analytic functions $\lambda_i(h)$ such that,
for any $h>0$
\begin{equation}
\hbox{spec}\ \Delta_h =\{\lambda_i(h):i=0,1,\ldots\}.
\end{equation}
The main result of the paper is an asymptotical formula
for the eigenvalue distribution function $N_h(\lambda)$
of the operator $\Delta_h$:
\begin{equation}
N_h(\lambda)=\sharp \{\lambda_i(h):\lambda_i(h)\leq \lambda\}.
\end{equation}
\begin{theorem}
\label{intr}
If $(M,{\cal F})$ be a Riemannian foliation,
equipped with a bundle-like Riemannian metric $g_M$.
Then the asymptotical formula for $N_h(\lambda)$
has the following form:
\begin{equation}
\label{eig1}
N_{h}(\lambda ) =h^{-q}
\frac{(4\pi)^{-q/2}}{\Gamma((q/2)+1)}
\int_{-\infty}^{\lambda} (\lambda - \tau )^{q/2}d_{\tau}N_{\cal F}(\tau
)+o(h^{-q}),
h\rightarrow 0,
\end{equation}
where $N_{\cal F}(\lambda )$ is the spectrum distribution function of the
tangential Laplace
operator
\begin{equation}
\Delta_F:C^{\infty}(M,\Lambda T^{*}M)
\rightarrow C^{\infty}(M,\Lambda T^{*}M).
\end{equation}
\end{theorem}
We refer the reader to Section~\ref{leaf} for a detailed
formulation of this Theorem.
We stated also the asymptotical formula for the trace
of an operator $f(\Delta_h)$ for any function
$f\in C_c(\Bbb R)$ (see Theorems~\ref{main} and \ref{fun}
below).

The study of asymptotical behaviour of geometric
objects (like as harmonic forms, eta-invariants etc.) associated with a family
of Riemannian
metrics on fibrations as the metrics become singular
was stimulated by Witten's work on adiabatic limits
\cite{W}. For further developments
see, for instance, \cite{Mel,Dai,F,Ge} and references there.

In the spectral theory of differential operators,
problems in question are related with the Born-Oppenheimer
approximation which consist in that
 the Schrodinger operator for polyatomic molecule is
 considered in the semiclassical limit where the mass ratio of electronic to
nuclear mass tends to zero (see, for instance,
\cite{KMS} and references there). In particular, the result
on semiclassical asymptotics for spectrum distribution function in a fibration
case is, essentially, due to
\cite{Bal}.

The investigation of semiclassical spectral asymptotics
for foliations
was started by the author in \cite{Ko1,Ko2,asymp}. There
we considered the problem in the operator setting, that is,
we studied spectral asymptotics for the self-adjoint
hypoelliptic operator $A_h$ of the form
\begin{equation}
A_h=A+h^mB,
\end{equation}
where $A$ is a tangentially elliptic operator of order
$\mu>0$
with the positive tangential principal symbol, and
$B$ be a differential operator of order $m$
on $M$ with the positive, holonomy invariant transversal
principal symbol and obtianed an asymptotical formula
for spectrum distribution function of this operator
when $h$ tends to zero.

In this work, we adapted our results on semiclassical
spectral asymptotics to the geometric setting
of adiabatic limits on foliations.

The main observation related with the asymptotical
formula (\ref{eig1}) is that its right-hand side
depends only on leafwise spectral data of the
tangential Laplace
operator $\Delta_F$. So, in a  case when the foliation
${\cal F}$ is nonamenable, there might to be a $\lambda>0$
such that
\begin{equation}
\label{eig2}
\lim_{h\rightarrow 0}h^qN_h(\lambda)=0.
\end{equation}
The formula~(\ref{eig2}) allows, in particular,
to introduce spectral
characteristics $r_k(\lambda)$ related with adiabatic limits which are
nontrivial in the nonamenable case.
We hope that some invariants of the function
$r_k(\lambda)$ introduced above  near $\lambda=0$ might to
be independent of the choice of
metric on $M$ (otherwise speaking, to be coarse invariants),
and, moreover, be topological or homotopic
invariants of foliated manifolds (just as in the case
of Novikov-Shubin invariants \cite{N-Sh}).
We discuss these questions and their relationships with the spectral theory of
leafwise
Laplacian and with noncommutative
spectral geometry of foliations in Section~\ref{disc}.

The organization of the paper is as follows.

In Section~\ref{pdo},
 we recall some facts on pseudodifferential operators
 on foliated manifolds.

In Section~\ref{geom}, we summarize some necessary
properties of the Laplace operator on a
foliated manifold.

In the Sections~\ref{funct} and \ref{proof},
we formulate and prove the asymptotical formula for
 $\hbox{tr} f(\Delta_h)$ when $h$ tends to
zero for any function $f\in C_c({\Bbb R})$.

In Section~\ref{leaf},
we rewrite the asymptotical formula of Section~\ref{funct}
in terms of spectral characteristics of the operator $\Delta_F$.
In particular, this provides a proof
of the main Theorem~\ref{intr} on an asymptotic behaviour of the
eigenvalue distribution function.

Finally, in Section~\ref{limits} we discuss some facts and
examples related with the asymptotical behaviour of individual
eigenvalues of the operator $\Delta_h$ when $h$ tends
to zero, and, as mentioned above, Section~\ref{disc}
is devoted to a discussion of various aspects of the main
asymptotical formula~(\ref{eig1}).

The work was done during a stay at the Max Planck Institut
f\"{u}r Mathematik at Bonn. I wish to express my
gratitude to it for hospitality and support.

\section{Pseudodifferential operators on foliations}
\label{pdo}
Here we recall some facts on pseudodifferential operators
 on foliated manifolds. The main
references here are \cite{tang,asymp}.

Let $(M,{\cal F})$ be a compact foliated manifold, $F$
be a distribution of tangent planes to ${\cal F}$.
The embedding $F\subset TM$ induces an embedding of
differential operators $Diff^\mu({\cal F})\subset
Diff^\mu(M)$, and differential operators on $M$ obtained
in such a way is called tangential differential operators.

More generally, let $E$ be an Hermitian vector bundle on $M$.
We say that a linear
differential operator $A$ of order $\mu$ acting  on  $C^{\infty}(M,E)$
is a tangential  operator, if, in
any foliated chart $\kappa\ :\ I^p\times I^q\rightarrow M$
($I=(0,1)$ is the open interval) and any
trivialization of the bundle $E$ over it, $A$ is of the form
\begin{equation}
A =\sum _{|\alpha|\leq \mu}a_{\alpha}(x,y)D^{\alpha}_{x}, (x,y)\in I^{p} \times
I^{q},
\label{(1.14)}
\end{equation}
with $a_{\alpha}$, being matrix valued functions on
$I^p\times I^q$.

Let $Diff^\mu({\cal F},E)$ denote the set of all tangential
differential operators of order $\mu$ acting in
$C^{\infty}(M,E)$.

Now we introduce the classes $Diff^{m,\mu}(M,{\cal F},E)$
by taking compositions of tangential differential
operators of order $\mu$ and differential operators of order
$m$ on $M$. That is, we say that $A\in Diff^{m,\mu}(M,{\cal F},E)$ if
$A$ is of the form
\begin{equation}
A=\sum_{\alpha} B_{\alpha}C_{\alpha},
\end{equation}
where $B_{\alpha}\in Diff^{m}(M,E)$, $C_{\alpha}\in
Diff^{\mu}({\cal F},E)$.

{}From symbolic calculus, it can be easily seen that:

(1) if
$A_1\in Diff^{m_1,\mu_1}(M,{\cal F},E)$,
$A_2\in Diff^{m_2,\mu_2}(M,{\cal F},E)$, then $A_1\circ
A_2\in Diff^{m_1+m_2,\mu_1+\mu_2}(M,{\cal F},E)$;

(2) if $A\in Diff^{m,\mu}(M,{\cal F},E)$, then the adjoint
$A^{*}\in Diff^{m,\mu}(M,{\cal F},E)$.

Classes $Diff^{m,\mu}(M,{\cal F},E)$ can be extended to
bigraded classes of pseudodifferential operators
$\Psi ^{m,\mu}(M,{\cal F},E)$, which contain, for instance,
parametrices for elliptic operators from the classes
$Diff^{m,\mu}(M,{\cal F},E)$.
We don't give its definition here, referring to
\cite{tang}(see also \cite{asymp}) for details
and will be
restricted by an introduction of classes of differential operators.

Now we recall the definition of a scale of Sobolev type spaces $H^{s,k}(M,{\cal
F},E), s\in
{\Bbb R}, k\in {\Bbb R}$, corresponding to classes of differential  operators
introduced  above.

The space $H^{s,k}({\Bbb R}^{n},{\Bbb R}^{p},{\Bbb C}^r)$
consists of all ${\Bbb C}^r$-valued tempered distributions
$u\in  S'({\Bbb R}^{n},{\Bbb C}^r)$ such that $\tilde{u} \in
L^{2}_{\hbox{loc}}({\Bbb R}^{n},{\Bbb C}^r)$ ($\tilde{u}$  the Fourier
transform) and
\begin{equation}
\label{(1.9)}
\| u\|^{2}_{s,k}
= \int \int  |\tilde{u}(\xi ,\eta )|^{2}(1 +|\xi |^{2} + |\eta |^{2})^{s} (1
+|\xi |^{2})^{k}d\xi  d\eta  < \infty.
\end{equation}
\noindent The identity (\ref{(1.9)}) serves as a definition of a norm $\|\  \|
_{s,k}$ in the
space $H^{s,k}({\Bbb R}^{n},{\Bbb R}^{p},{\Bbb C}^r)$.

The space $H^{s,k}(M,{\cal F},E)$ consists of all $u\in  {\cal D}'(M,E)$ such
that,  for
any foliated coordinate chart $\kappa  : I^{p} \times I^{q} \rightarrow  U =
\kappa (I^{p} \times I^{q}) \subset M$,
any trivialization of the bundle $E$ over it, and
for any $\phi  \in C^{\infty }_{c}(U)$,  the  function $\kappa^{\ast}(\phi u)$
belongs  to  the  space
$H^{s,k}({\Bbb R}^{n},{\Bbb R}^{p},{\Bbb C}^r)$ ($r=\hbox{rank}\ E$). Fix some
finite covering $\{ U_{i} : i = 1,\ldots ,d \}$ of $M$  by
foliated coordinate patches with the foliated coordinate charts $\kappa _{i} :
I^{p} \times I^{q} \rightarrow U_{i} = \kappa _{i}(I^{p} \times I^{q})$ and
trivializations of the bundle
$E$ over them, and a partition of unity $\{ \phi _{i} \in C^{\infty }(M):  i
= 1,\ldots,d \}$ subordinate to  this  covering.  A  scalar  product  in
$H^{s,k}(M,{\cal F},E)$ is defined by the formula
\begin{equation}
(u,v)_{s,k} =\sum ^{d}_{i=1} (\kappa^{\ast}(\phi _{i}u), \kappa ^{\ast}(\phi
_{i}v))_{s,k}, u,v\in H^{s,k}(M,{\cal F},E).
\label{(1.10)}
\end{equation}
We have the  following result on the action of   differential
operators  of  class $Diff^{m,\mu}(M,{\cal F},E)$ in the spaces
$H^{s,k}(M,{\cal F},E)$ (see \cite{tang,asymp} for a proof
in the scalar case).

\begin{proposition}
\label{Proposition 1.3}
An operator $A\in Diff ^{m,\mu}(M,{\cal F},E)$ defines  a  linear
bounded operator from $H^{s,k}(M,{\cal F},E)$ to $H^{s-m,k-\mu}(M,{\cal F},E)$
for any $s\in {\Bbb R}$, $k\in
{\Bbb R}$.
\end{proposition}

Finally, the scale of Sobolev type spaces introduced above
allows us to formulate a Garding inequality for tangentially
elliptic operators (for  the  proof, see \cite{tang}).

\begin{proposition}
\label{tangG}
If $A$ is  tangentially  elliptic  operator  of
order $\mu $ with the positive tangential principal  symbol,  then,  for
any $s\in {\Bbb R}, k\in {\Bbb R}$, there exist constants $C_{1}>0$ and $C_{2}$
such that
\begin{equation}
\hbox{Re}\ (Au,u)_{s,k} \geq C_{1}\| u\| ^{2}_{s,k+\mu/2} - C_{2}\| u\|
^{2}_{s,-\infty}, u\in  C^{\infty}(M,E).
\label{(1.18)}
\end{equation}
\end{proposition}

\section{Geometric operators on Riemannian foliations}
\label{geom}
Here we summarize some necessary
properties of the Laplace operator on a
foliated manifold.

As above, $(M,{\cal F})$ denotes a closed foliated Riemannian manifold,
$\dim M = n$, $\dim {\cal F} = p$, $p+q=n$,
equipped with a Riemannian metric $g_M$, $F=T{\cal F}$
be an integrable distribution
of $p$-planes in $TM$. Recall that we choose the orthogonal
complement  $H$ to $F$, so
\begin{equation}
\label{decomp1}
F\bigoplus H=TM.
\end{equation}
The decomposition~(\ref{decomp1}) induces a bigrading on $\Lambda T^{*}M$
by the formula
\begin{equation}
\label{bigrad}
\Lambda^k T^{*}M=\bigoplus_{i=0}^{k}\Lambda^{i,k-i}T^{*}M,
\end{equation}
where
\begin{equation}
\Lambda^{i,j}T^{*}M=\Lambda^{i}F^{*}\bigotimes
\Lambda^{j}H^{*}.
\end{equation}

Now we transfer the family $\Delta_h$ to a fixed Hilbert space $L^{2}(M,\Lambda
T^{*}M ,g)$. For this goal we introduce the isometry
\begin{equation}
\Theta_{h} : L^{2}(M,\Lambda T^{*}M ,
g_{h})\rightarrow  L^{2}(M,\Lambda T^{*}M ,g),
\end{equation}
where, for $u \in  L^{2}(M,\Lambda^{i,j}T^{*}M , g_{h})$, we have
\begin{equation}
\Theta_{h}u = h^{j}u.
\end{equation}
The operator $\Delta_h$ in the Hilbert space
$L^{2}(M,\Lambda T^{*}M , g_{h})$ corresponds under the isometry
$\Theta_{h}$ to the operator
\begin{equation}
L_{h}=
\Theta_{h}\Delta_h\Theta_{h}^{-1}
\end{equation}
in the Hilbert space
$L^2(M,\Lambda T^{*}M)=L^{2}(M,\Lambda T^{*}M ,g)$.

De Rham differential $d$ inherits the
decomposition~(\ref{decomp1}) in the form
\begin{equation}
d=d_F+d_H+\theta.
\end{equation}
Here the tangential de Rham differential $d_F$
and the transversal de Rham differential $d_H$ are first order differential
operators, and $\theta$ is zeroth order. Moreover, the
operator $d_F$ doesn't depend on a choice of the orthogonal
complement $H$ (see, for instance, \cite{Re}).

Then we have the following assertion on the form of the
operator $L_h$.

\begin{lemma}[\cite{F}]
\label{Lh}
We have
\begin{equation}
L_{h}=d_h\delta_h + \delta_h d_h,
\end{equation}
where
\begin{equation}
d_h = d_F + hd_H + h^{2}\theta,
\end{equation}
and
\begin{equation}
\delta_h = \delta_F + h \delta_H + h^{2}\theta^{*},
\end{equation}
is the adjoint, where $\delta_F$, $ \delta_H$ and $ \theta^{*}$
are the adjoints to
$d_F$, $d_H$ and $\theta$ respectively. Here
we consider the adjoints taken in the Hilbert space
$L^2(M,\Lambda T^{*}M)$.
\end{lemma}

By Lemma~\ref{Lh}, the operator $L_h$ is of the
following form:
\begin{equation}
\label{L}
L_h =\Delta_F + h^2\Delta_H + h^4\Delta_{-1,2}+ hK_1+h^2K_2 +h^3K_3,
\end{equation}
where
\begin{itemize}
\item The operator
\begin{equation}
\Delta_F=d_F\delta_F+\delta_Fd_F\in Diff^{0,2}(M,{\cal F},
\Lambda T^{*}M )
\end{equation}
is the tangential
Laplacian in the space $C^{\infty}(M,\Lambda T^{*}M)$.
\item The operator
\begin{equation}
\Delta_H=d_H \delta_H+ \delta_Hd_H\in Diff^{2,0}(M,{\cal F},\Lambda T^{*}M )
\end{equation}
is the transversal Laplacian in the space $C^{\infty}(M,\Lambda T^{*}M)$.
\item $\Delta_{-1,2}=\theta\theta^{*}+ \theta^{*}\theta
\in Diff^{0,0}(M,{\cal F}, \Lambda T^{*}M ))$.
\item $K_1 =  d_F \delta_H+ \delta_H
d_F + \delta_Fd_H+d_H
\delta_F\in Diff^{1,0}(M,{\cal F}, \Lambda T^{*}M ))$.
\item $K_2 =  d_F \theta^{*}+ \theta^{*}
d_F + \delta_F\theta+\theta
\delta_F\in Diff^{0,0}(M,{\cal F}, \Lambda T^{*}M )).$
\item $K_3 = d_H \theta^{*}+ \theta^{*}
d_H +  \delta_H\theta+\theta
 \delta_H\in Diff^{1,0}(M,{\cal F}, \Lambda T^{*}M ))$.
\end{itemize}

{}From now on, we will assume that $(M,{\cal F})$ is
a Riemannian foliation with a bundle-like metric $g_{M}$, that
is, it satisfies one of the following equivalent
conditions (see \cite{Re}):
\begin{enumerate}
\item $(M,{\cal F})$ locally has the structure of Riemannian
submersion;
\item for any $X\in F$ we have
\begin{equation}
\label{inv}
\nabla_X^{\cal F}g_H=0,
\end{equation}
where $\nabla^{\cal F}$ is a Bott connection on $H$;
\item the distribution $H$ is totally geodesic.
\end{enumerate}

\noindent The following Lemma states the main specific property
of geometrical operators on Riemannian foliated manifold.

\begin{lemma}
\label{op}
If $(M,{\cal F})$ is
a Riemannian foliation with a bundle-like metric $g_{M}$,
then the operators
\begin{equation}
d_F \delta_H+ \delta_Hd_F\ \hbox{and}\  \delta_Fd_H+d_H\delta_F
\end{equation}
belong to the class $Diff^{0,1}(M,{\cal F}, \Lambda T^{*}M ))$. In particular,
we have
\begin{equation}
K_1\in Diff^{0,1}(M,{\cal F}, \Lambda T^{*}M )).
\end{equation}
\end{lemma}

For any $h>0$, the operator $L_h$ is a formally self-adjoint, elliptic
operator in $L^2(M,\Lambda T^{*}M)$ with the positive principal symbol. The
following Proposition
is a refinement of the classical Garding inequality
for the operator $L_h$ in $H^{s,k}(M,{\cal F},\Lambda T^{*}M )$
\begin{proposition}
\label{Garding}
Under current hypotheses, there exists constants
$C_1>0$, $C_2>0$ and $C_3>0$ such that for any $h>0$
small enough we have the following inequality:
\begin{equation}
\label{G1}
(L_hu,u)\geq (1-C_1h^2)(\Delta_Fu,u)+C_2h^2\|u\|^2_{1,0}
-C_3\|u\|^2, u\in C^{\infty}(M,\Lambda T^{*}M )).
\end{equation}
\end{proposition}

\noindent{\bf Proof.}\
By (\ref{L}), we have
\begin{eqnarray}
(L_hu.u) &=&(\Delta_Fu,u) + h^2(\Delta_Hu,u) + h^4(\Delta_{-1,2}u,u)\nonumber\\
&+& h(K_1u,u)+h^2(K_2u,u) +h^3(K_3u,u),\nonumber\\
& & u\in C^{\infty}(M,\Lambda T^{*}M).
\end{eqnarray}
It clear that $(\Delta_{-1,2}u,u)\geq 0$.
By Proposition~\ref{Proposition 1.3}, we have
\begin{eqnarray}
\label{K2}
(K_2u,u)\geq -C_4\|u\|^2,
(K_3u,u)\geq -C_5\|u\|^2_{1,0}.
\end{eqnarray}
So we obtain
\begin{eqnarray}
\label{Lh2}
(L_hu.u) & \geq &(\Delta_Fu,u) + h^2(\Delta_Hu,u)
+ h(K_1u,u)\nonumber\\
&-& C_4h^2\|u\|^2 -C_5h^3\|u\|^2_{1,0}.
\end{eqnarray}
The operator $\Delta_F+\Delta_H$ is an second order
elliptic operator with the positive principal symbol,
so, by the standard Garding inequality, we have
\begin{equation}
\label{01}
((\Delta_F+\Delta_H)u,u)\geq C_6\|u\|^2_{1,0}-
C_7\|u\|^2,
\end{equation}
that implies the estimate
\begin{eqnarray}
\label{Lh3}
(L_hu.u) & \geq & (1-h^2)(\Delta_Fu,u) +C_7h^2\|u\|^2_{1,0}+ h(K_1u,u)
\nonumber\\
&-& C_8\|u\|^2.
\end{eqnarray}

Finally, we make use of the inequality
\begin{equation}
\label{K1}
|(K_1u,u)| \leq  C_7\|u\|_{0,1}\|u\| \leq C_8(h\|u\|^2_{0,1}+h^{-1}\|u\|^2)
\end{equation}
and the tangential Garding estimate (see Proposition
\ref{tangG})
\begin{equation}
\|u\|^2_{0,1}\leq C_9((\Delta_Fu,u)+ \|u\|^2),
\end{equation}
that completes immediately the proof.
\medskip
\par
\noindent{\bf Remark.}\  In some cases, it is sufficient
 to use more crude estimate
\begin{equation}
\label{crude}
(L_hu,u)\geq C_1\|u\|^2_{0,1}+C_2h^2\|u\|^2_{1,0}
-C_3\|u\|^2, u\in C^{\infty}(M,\Lambda T^{*}M )),
\end{equation}
which follows from (\ref{G1}), if we apply the standard
Sobolev norm estimate
\begin{equation}
(\Delta_Fu,u)\leq C_{10}\|u\|^2_{0,1}.
\end{equation}

Let $H_h(t)=\exp (-tL_h),t\geq 0,$ be the parabolic semigroup, generated by the
operator $L_h$:
\begin{equation}
H_h(t):C^{\infty}(M,\Lambda T^{*}M)
\rightarrow C^{\infty}(M,\Lambda T^{*}M).
\end{equation}
 For any $t>0$, the operator $H_h(t)$ is an operator with a smooth kernel.
Proposition~\ref{Garding} implies
 the following norm estimates for operators of
this semigroup in the spaces
$H^{s,k}(M,{\cal F},\Lambda T^{*}M )$
(see also \cite{asymp}).
\begin{proposition}
\label{par}
We have the following estimates:
\begin{equation}
\|H_h(t)u\|  _{r,k} \leq C_{rsk}t^{(s-k-r)/2}h^{s-r}\| u\| _{s} , u\in
C^{\infty}(M,\Lambda T^{*}M ),
\end{equation}
\noindent if \(r>s, h\in (0,1], 0< t\leq 1\),
and the estimate
\begin{equation}
\|H_h(t)u\| _{s,k} \leq  C_{sk} t^{-k/2}\| u\| _{s}, u\in C^{\infty}(M,\Lambda
T^{*}M ).
\end{equation}
\noindent if \(r = s, h\in [0, 1], 0< t\leq 1\), where the constants don't
depend on $t$ and $h$.
\end{proposition}

\section{Asymptotical formula for the functions
of the Laplace operator}
\label{funct}
Form now on, we will assume that $(M,{\cal F})$ is a  Riemannian foliation,
equipped with a bundle-like Riemannian metric $g_M$.
In this Section, we state the asymptotical formula for
 $\hbox{tr} f(\Delta_h)$ when $h$ tends to
zero for any function $f\in C_c({\Bbb R})$.

We will denote by $G_{\cal F}$ the holonomy groupoid of $(M,{\cal F})$. Recall
that $G_{\cal F}$ is equipped with the source and the target
maps $s,r:G_{\cal F}\rightarrow M$. We will make use of the standard
notation:
$G_{\cal F}^{(0)}=M$ is the set of objects, $G_{\cal F}^x=\{\gamma\in G_{\cal
F}:
r(\gamma)=x\}$,$x\in M$. Recall that $G_{\cal F}^x$ is the covering
of the leaf through the point $x$, associated with the
holonomy group of the leaf. We will identify a point $x\in  M$ with the
identity element in $G^{x}_{x}$.
Finally, we will denote by $\lambda _{L}$ the Riemannain volume form on each
leaf $L$ of ${\cal F}$ and by $\lambda ^{x}$
its lift to a measure
on the holonomy covering  $G_{\cal F}^x, x\in  M$.

For any vector bundle $E$ on $M$, we denote by
$C^{\infty}_c(G_{\cal F},E)$ the space of all smooth, compactly supported
sections of the vector bundle
$(s,r)^*(E^*\bigotimes E)$ over $G_{\cal F}$. In other
words, for any $k\in C^{\infty}_c(G_{\cal F},E)$, its
value at a point $\gamma\in G_{\cal F}$ is a linear
map $k(\gamma):E_{s(\gamma)}\rightarrow E_{r(\gamma)}$.
We will use a correspondence between tangential kernels
$k\in C^{\infty}_c(G_{\cal F},E)$ and tangential operators
$K:C^{\infty}(M,E)\rightarrow C^{\infty}(M,E)$ via the
formula
\begin{equation}
Ku(x)=\int_{G_{\cal F}^x}k(\gamma)u(s(\gamma))d\lambda^x(\gamma),
u\in C^{\infty}(M,E).
\end{equation}

Now we introduce
a notion of a principal $h$-symbol of the operator
$\Delta_h$. It is well-known (see, for instance,
\cite{Molino,Re}) that the conormal bundle
$H^{*}$ to the foliation ${\cal F}$ has a partial (Bott)
connection, which is flat along the leaves of the foliation.
So we can lift the foliation ${\cal F}$ to the
foliation ${\cal F}_H$ in the conormal bundle $H^{*}$. The leaf
$\tilde{L}_{\nu}$ of the foliation ${\cal F}_H$  through a point $\nu\in H^{*}$
is diffeomorphic to the holonomy
covering $G_{\cal F}^x$ of the leaf $L_x, x=\pi(\nu)$
of the foliation ${\cal F}$ through the point $x$ (here $\pi:H^{*}\rightarrow
M$ is the bundle map) and has a trivial holonomy.

Denote by
\begin{equation}
\Delta_{{\cal F}_H}:C^{\infty}(H^{*},\pi^{*}
\Lambda T^{*}M)\rightarrow C^{\infty}(H^{*},\pi^{*}
\Lambda T^{*}M)
\end{equation}
the lift of the
leafwise Laplacian $\Delta_F$ to tangentially
elliptic operator on $H^{*}$ with respect to
${\cal F}_H$.
\medskip
\par
\noindent{\bf Remark.} If we fix $x\in M$, the restriction of the foliation
${\cal F}_H$ on $H^{*}_x$ is a linear model of the foliation
${\cal F}$ in some neighborhood of the leaf $L_x$
through a point
$x$, so the restriction $\Delta_x$ of $\Delta_{{\cal F}_H}$
on $H^{*}$,
\begin{equation}
\Delta_x:C^{\infty}(H^{*}_x,\pi^{*}
\Lambda T^{*}L\bigotimes \Lambda H^{*}_x)\rightarrow C^{\infty}(H^{*}_x,\pi^{*}
\Lambda T^{*}L\bigotimes \Lambda H^{*}_x),
\end{equation}
is the model operator for the tangential Laplacian
$\Delta_F$ at the "point" $x\in M/{\cal F}$.
\medskip
\par
\noindent {\bf Definition.} The principal $h$-symbol  of the operator
$\Delta_h$ is a tangentially elliptic
operator
\begin{equation}
\sigma_h(\Delta_h):C^{\infty}(H^{*},\pi^{*}
\Lambda T^{*}M)\rightarrow C^{\infty}(H^{*},\pi^{*}
\Lambda T^{*}M)
\end{equation}
on $H^{*}$ with respect to the foliation
${\cal F}_H$, given by the formula
\begin{equation}
\sigma_h(\Delta_h) = \Delta_{{\cal F}_H}+g_H,
\end{equation}
where $g_H$ is the scalar multiplication operator by
the function $g_H(\nu),\nu \in H^{*}$.
\medskip
\par
The holonomy groupoid $G_{{\cal F}_H}$ of the lifted
foliation ${\cal F}_H$ consists of all triples
$(\gamma,\nu,\eta)\in G_{\cal F}\times H^{*}\times H^{*}$ such
that $s(\gamma)=\pi(\nu)$, $r(\gamma)=\pi(\eta)$ and
$(dh_{\gamma}^{*})^{-1}(\nu)=\eta$, where
$dh_{\gamma}^{*}$ is codifferential of the holonomy
map, with the source map $s:G_{{\cal F}_H}\rightarrow
H^{*}, s(\gamma,\nu,\eta)=\nu$ and the target map
$r:G_{{\cal F}_H}\rightarrow
H^{*}, r(\gamma,\nu,\eta)=\eta$.The projection
$\pi:H^{*}\rightarrow M$ induces the map
$\pi_G:G_{{\cal F}_H}\rightarrow G_{\cal F}$ by
\begin{equation}
\pi_G(\gamma,\nu,\eta)=\gamma, (\gamma,\nu,\eta)\in
G_{{\cal F}_H}.
\end{equation}
Denote by $\hbox{tr}_{{\cal F}_H}$ the trace on
the von Neumann algebra $W^{\ast}(G_{{\cal F}_H},\pi^{*}\Lambda T^{*}M)$
of all tangential operators on $H^{*}$ with respect to the foliation ${\cal
F}_H$,
given by a holonomy invariant measure $dx\ d\nu$ on
$H^{*}$ \cite{Co3}. For any tangentially elliptic operator $K$
on $(H^{*},{\cal F}_H)$, given by the tangential kernel
$k\in C^{\infty }_c(G_{{\cal F}_H},\pi^{*}\Lambda T^{*}M)$,
$k=k(\gamma,\nu,\eta)$
we have
\begin{equation}
\hbox{tr}_{{\cal F}_H}(K)=\int_{H^{\ast}}
 \hbox{Tr}_{\pi^{*}\Lambda T^{*}M}\ k(x,\nu,\nu)dxd\nu.
\end{equation}

\begin{theorem}
\label{main}
For  any function $f\in C_c({\Bbb R})$, we  have the asymptotical formula
\begin{equation}
\label{f}
\hbox{tr}\ f(\Delta_{h}) =(2\pi)^{-q}h^{-q} \hbox{tr}_{{\cal F}_H}\
f(\sigma_h(\Delta_h))
+O(h^{1-q}), h\rightarrow 0.
\end{equation}
\end{theorem}
We will prove this Theorem in the next Section, and now
we conclude the Section with some remarks.
\medskip
\par
\noindent{\bf Remarks.}\ (1)\ In a case of the Schrodinger
operator on a compact manifold $M$ with an operator-valued potential $V\in
{\cal L}(H)$ with a Hilbert space $H$ such
that
$V(x)^{*}=V(x)$ (a fibration
case)
\begin{equation}
H_h=-h^2\Delta +V(x), x\in M,
\end{equation}
the corresponding asymptotical formula has the following form:
\begin{equation}
\label{semi1}
\hbox{tr}\ f(\Delta_h)=(2\pi)^{-n}h^{-n}\int\ Tr\
f(h(x,\xi))dxd\xi+o(h^{-n}), h\rightarrow 0+,
\end{equation}
where $h(x,\xi)$ is the operator-valued principal
$h$-symbol
\begin{equation}
h(x,\xi)=|\xi|^2+V(x), (x,\xi)\in T^{*}M.
\end{equation}
So the formula~(\ref{f}) has the same form as (\ref{semi1})
with the difference that the usual integration over the base
and the fibrewise trace are replaced by the integration
in a sense of the noncommutative integration theory
\cite{Co3}.
\medskip
\par
\noindent (2)\ We don't make an essential use of
 a operator-valued symbolic calculus. Indeed, it
is a difficult problem to develop such a calculus
in a general case. The introduction of the principal
$h$-symbol of the operator $\Delta_h$ allow us to
simplify the final asymptotical formula and also
some algebraic calculations (see below for a passage
from an asymptotical formula for $\hbox{tr}\ \exp(-t\Delta_h)$ to an
asymptotical formula for $\hbox{tr}\ f(\Delta_h)$ with an arbitrary function
$f\in C^{\infty}_c({\Bbb R})$.

\section{Proof of Theorem~\protect{\ref{main}}}
\label{proof}
In this Section, we prove Theorem~\ref{main}, concerning
an asymptotical behaviour of $\hbox{tr}\ f(\Delta_h)$ when $h$ tends
to zero.

First of all,
let us note that, without loss of generality, we may consider an asymptotical
behaviour
of $\hbox{tr}\ f(L_h)$.
The proof of Theorem~\ref{main} relies on a comparison of the operator $L_h$
with
some operator $\bar{L}_h$  of the almost product structure as in \cite{asymp}
with a subsequent use of
results of \cite{asymp} on semiclassical
spectral asymptotics for elliptic operators on foliated manifolds.

So let the operator $\bar{L}_h\in
Diff^{2,0}(M,{\cal F}, \Lambda T^{*}M ))$
be given by the formula
\begin{equation}
\bar{L}_h =\Delta_F + h^2\Delta_H.
\end{equation}
The operators $L_h$ and $\bar{L}_h$ are generators
parabolic semigroups of linear bounded operators
in the space $L^2(M,\Lambda T^{*}M)$ denoted by
\begin{eqnarray}
H_h(t)&=&e^{-tL_h}, t\geq 0,\\
\bar{H}_h(t)&=&e^{-t\bar{L}_h}, t\geq 0,
\end{eqnarray}
respectively. It is clear that,
indeed, these operators are  smoothing operators when
$t>0$.

The operator $\bar{L}_h$ satisfies the conditions
of \cite{asymp}, that is, it is of the form
\begin{equation}
\bar{L}_h=A+h^2B,
\end{equation}
where $A=\Delta_F$ is a second order tangentially elliptic operator
with the scalar, positive tangential principal symbol, and
$B=\Delta_H$ be a second order differential operator
on $M$ with the scalar, positive, holonomy invariant transversal
principal symbol. Indeed, it is easy to see that the transversal principal
symbol of operator $\Delta_H$,
which is the restriction of its principal symbol
from $T^{*}M$ to
the conormal bundle $H^{*}$, is given by the formula
\begin{equation}
\sigma (\nu) =g_{H^{*}}(\nu)I, \nu\in H^{*},
\end{equation}
and its holonomy invariance is equivalent to the assumption
on the foliation ${\cal F}$ to be Riemannian (see (\ref{inv})).
\medskip
\par
\noindent{\bf Remark.}
The only necessary property which we need from holonomy
invariance condititon is the fact that the commutator
$[A,B]$, which, by general symbolic calculus, belongs
to the class $Diff^{2,1}(M,{\cal F}, \Lambda T^{*}M )$, is an operator of the
class
$Diff^{1,2}(M,{\cal F}, \Lambda T^{*}M )$,
and this fact can be checked by a straightforward
calculation and looks very similar to the second
assertion of Lemma~\ref{op}.
\medskip
\par
By \cite{asymp}, the operators of the parabolic semigroup $\bar{H}_h(t)$
satisfy the same estimate
as in Proposition~\ref{par}.
\begin{equation}
\label{par1}
\| \bar{H}_h(t)u\| _{r,k} \leq C_{r,s,k}t^{(s-k-r)/2}h^{s-r}\| u\| _{s} , u\in
C^{\infty}(M,\Lambda T^{*}M ),
\end{equation}
\noindent if \(r>s, h\in (0,1], 0< t\leq 1\), and the estimate
\begin{equation}
\| \bar{H}_h(t)u\| _{s,k} \leq  C_{sk} t^{-k/2}\| u\| _{s}, u\in
C^{\infty}(M,\Lambda T^{*}M ).
\label{(3.13)}
\end{equation}
\noindent if \(r = s, h\in [0, 1], 0< t\leq 1\), where the constants
don't depend on $t$ and $h$.

Now we want to compare the semigroups  $\bar{H}_h(t)$
and $H_h(t)$. First, we state the norm estimates for the difference $H_h(t) -
\bar{H}_h(t)$.

\begin{proposition}
\label{dif}
We have the estimate
\begin{equation}
\label{dl00}
\|(H_h(t) -  \bar{H}_h(t))u\| _{r,k} \leq C_{r,s,k}t^{(s-k-r)/2}h^{s-r-1}\| u\|
_{s} , u\in C^{\infty}(M,\Lambda T^{*}M ),
\end{equation}
\noindent if \(r>s, h\in (0,1], 0< t\leq 1\),
and the estimate
\begin{equation}
\label{dl01}
\|(H_h(t) -  \bar{H}_h(t))u\| _{s,k} \leq  C_{sk} t^{-k/2}\| u\| _{s}, u\in
C^{\infty}(M,\Lambda T^{*}M ).
\end{equation}
\noindent if \(r = s, h\in [0, 1], 0< t\leq 1\), where the constants
don't depend on $t$ and $h$.
\end{proposition}
\noindent {\bf Proof.} For a proof, we make use of
the Duhamel formula
\begin{equation}
(H_h(t) -  \bar{H}_h(t))u =
\label{(4.4)}
\int^{t}_{0} H_h(\tau)( \bar{L}_h-  L_h)
\bar{H}_h(t-\tau)u d\tau .
\end{equation}
We know the norm estimates for operators $H_h(t)$ and $\bar{H}_h(t)$ (see
Propositions~\ref{par} and (\ref{par1}))
and the explicit formula for the difference
$\bar{L}_h-  L_h$:
\begin{equation}
\label{dl}
L_h-\bar{L}_h =h^4\Delta_{-1,2}+ hK_1+h^2K_2 +h^3K_3.
\end{equation}
from where Proposition is proved in a usual way.
\medskip
\par
Now we pass from the Sobolev estimates for the
operator $H_h(t) -  \bar{H}_h(t)$ to pointwise and
trace estimates.
\begin{proposition}
\label{diftr}
Under current hypotheses, we have the estimates
\begin{equation}
|\hbox{tr}(H_h(t)-\bar{H}_h(t))|\leq Ch^{1-q}.
\end{equation}
\end{proposition}
\noindent{\bf Proof.}
For the proof,
we make use the following proposition
(see \cite{asymp} for a scalar case):
\begin{proposition}
\label{Proposition 1.5}
Let $(M,{\cal F})$ be a compact foliated manifold,
$E$ be an Hermitian vector bundle on $M$.
For any $s>p/2$ and  $k>q/2$, there  is  a
continuous embedding
\begin{equation}
H^{s,k}(M,{\cal F},E )
\subset C(M,E).
\label{(1.12)}
\end{equation}
\noindent Moreover, for any $s> p/2$ and $k> q/2$, there  is  a  constant
$C_{s,k} > 0$ such that, for each $\lambda  \geq  1$,
\begin{equation}
\label{(1.13)}
\sup_{x\in M} \vert u(x)\vert
\leq C_{s,l}\lambda ^{q/2}(\lambda ^{-s}\| u\| _{s,k} + \| u\| _{0,k+s}), u\in
H^{s,k}(M,{\cal F},E).
\end{equation}
\end{proposition}

Denote by $H_h(t,x,y)$ ($\bar{H}_h(t,x,y)$) the
integral kernels of operators $H_h(t)$ ($\bar{H}_h(t)$)
respectively.
Then, by
Propositions~\ref{dif} and \ref{Proposition 1.5}, we  obtain:
\begin{equation}
|H_h(t,x,x) -  \bar{H}_h(t,x,x)| \leq Ch^{1-q},
x\in M.
\end{equation}
that, due to the well-known formula for the trace
of an integral operator $K$ in the Hilbert space
$L^2(M,\Lambda T^{*}M)$ with a smooth kernel $k(x,y)$:
\begin{equation}
\hbox{tr}\ K =\int_M \hbox{Tr}\ k(x,x)dx,
\end{equation}
immediately completes the proof.
\medskip
\par
Denote by $h_{\cal F}(t,\gamma) \in
C^{\infty}(G_{\cal F},\Lambda T^{*}M)$ the tangential kernel of  the  smoothing
tangential operator $\exp(-t\Delta_F)$.

\begin{proposition}
\label{Theorem 4.1}
For  any $t>0$, we  have the asymptotical formula
\begin{equation}
\label{(4.2)}
\hbox{tr}\ e^{-tL_{h}} =(2\pi)^{-q}
h^{-q} \int_{M}(\int_{H^{\ast}_x} e^{-tg_H(\nu)}  d\nu )\ \hbox{Tr}_{\Lambda
T^{*}M}\ h_{\cal F}(t,x)dx+O(h^{1-q}), h\rightarrow 0.
\end{equation}
\end{proposition}
\noindent{\bf Proof.} By Propositions~\ref{par} and \ref{Proposition 1.5},
we have the estimate
\begin{equation}
\hbox{tr}\ e^{-tL_{h}}\leq Ch^{-q}, h\rightarrow 0.
\end{equation}
Moreover, by Proposition~\ref{diftr}, asymptotics of traces of  the operators
$H_h(t)$
and $\bar{H}_h(t)$  when $h$ tends to zero have the same leading terms
(of order $h^{-q}$), and we can apply the asymptotical
formula of \cite{asymp}
to complete the proof.
\medskip
\par
\noindent{\bf Remarks.} (1)\ Since
\begin{equation}
\int_{H^{\ast}_x} e^{-tg_H(\nu)}d\nu={\pi}^{q/2}t^{-q/2},
\end{equation}
the formula~(\ref{(4.2)}) can be rewritten in a
simpler form:
\begin{equation}
\label{simple}
\hbox{tr}\ e^{-tL_{h}} =(4\pi t)^{-q/2}
h^{-q} \int_{M}\ \hbox{Tr}_{\Lambda T^{*}M}\ h_{\cal F}(t,x)dx+O(h^{1-q}),
h\rightarrow 0.
\end{equation}
{}From~(\ref{simple}), we can also obtain an asymptotical
formula for the spectrum distribution function, but
it is more convenient for us to use the formula
in the form~(\ref{(4.2)}).
\medskip
\par
\noindent (2)\ For any $x\in M$,
the restriction $h_{\cal F}(t,\gamma) \in
C^{\infty}(G^x_{\cal F},\Lambda T^{*}M)$ of $h_{\cal F}$
on $G^x_{\cal F}$ is the kernel of the operator $\exp(-t\Delta_x)$, where
$\Delta_x$ the restriction of
$\Delta_F$ on $G^x_{\cal F}$ (see also Section~\ref{leaf}).
This fact doesn't extend to more general functions
$f(\Delta_F)$ (see \cite{tang}), and this is closely
related with so-called spectrum coincidence theorems
and with appearance of nonstandard asymptotical
formula~(\ref{eig2}).
\medskip
\par
\noindent{\bf Proof of Theorem~\ref{main}.}
The tangential kernel $h_{{\cal F}_H}(t)
\in C^{\infty}(G_{{\cal F}_H},\pi^{*}\Lambda T^{*}M)$ of the operator
$\exp(-t\Delta_{{\cal F}_H})$
is related with the tangential kernel $h_{{\cal F}}(t)\in C^{\infty}(G_{\cal
F},\Lambda T^{*}M)$
of operator
$\exp (-t\Delta_F)$ by the formula
\begin{equation}
h_{{\cal F}_H}(t,\gamma,\nu,\eta)=\pi_G^{*}h_{\cal F}(t,\gamma).
\end{equation}
The essential difference of the case of Riemannian
foliation from the general one consists in the
fact that the operators $\Delta_{{\cal F}_H}$ and $g_H$
considered as operators on $H^{*}$ commutes. In
particular, we have
\begin{equation}
e^{-t\sigma_h(\Delta_h)}= e^{-tg_H(\nu)}e^{-t\Delta_{{\cal F}_H}}, t>0.
\end{equation}
So the formula (\ref{(4.2)}) can be rewritten in terms of
the notation of this Section  as follows:
\begin{equation}
\hbox{tr}\ e^{-tL_{h}} =
h^{-q} \hbox{tr}_{{\cal F}_H}\ e^{-t\sigma_h(\Delta_h)}+O(h^{1-q}),
h\rightarrow 0.
\end{equation}
{}From where, using standard approximation arguments,
the theorem follows immediately.
\medskip
\par
\noindent{\bf Remark.}\ The passage from the operator $L_h$ to the operator
$\bar{L}_h$ resembles the passage from the Riemannian
connection on $M$ to the almost product connection
as in \cite{AL-T,Re}.

\section{Formulation in terms of leafwise spectral
characteristics}
\label{leaf}
Here we will write the asymptotical formula~(\ref{f})
in terms of spectral characteristics of the operator $\Delta_F$.
In particular, we obtain a proof
of the main theorem on an asymptotic behaviour of the
eigenvalue distribution function.

Recall that $\Delta_F$ denotes the tangential Laplacian
in the space $C^{\infty}(M,\Lambda T^{*}M)$.
Let us restrict the operator
$\Delta_F$ to the leaves of the foliation ${\cal F}$ and lift the restrictions
to
holonomy coverings of leaves. We obtain the family
\begin{equation}
\Delta_x : C^{\infty}_c(G_{\cal F}^x,r^{*}\Lambda T^{*}M)\rightarrow
C^{\infty}_c(G_{\cal F}^x,r^{*}\Lambda T^{*}M)
\end{equation}
of Laplacians on holonomy coverings of leaves.
By the hypotheses of Riemannian foliation,  the operator $\Delta_x$ is formally
self-adjoint in
$L^2(G_{\cal F}^x,r^{*}\Lambda T^{*}M)$, that, in turn, implies its
essential self-adjointness in this Hilbert space (with  initial  domain
$C^{\infty}_c(G_{\cal F}^x,r^{*}\Lambda T^{*}M)$) for any $x\in M$.
For each $\lambda  \in {\Bbb R}$, the kernel $e(\gamma,\lambda ), \gamma\in
G_{\cal F}$ of the spectral projections
 of the operators
$\Delta_x$, corresponding to the semiaxis $(-\infty ,\lambda ]$ define an
element of the von Neumann algebra
$W^{*}(G_{\cal F},\Lambda T^{*}M)$.  The section $e(\gamma,\lambda )$ is  a
leafwise smooth section of the bundle $(s^{*}\Lambda T^{*}M)^{*}
\bigotimes r^{*}\Lambda T^{*}M$ over
$G_{\cal F}$.

We introduce the spectrum distribution function
$N_{\cal F}(\lambda)$ of the operator $\Delta_F$ by the formula
\begin{equation}
N_{\cal F}(\lambda)=\int_M \hbox{Tr}_{\Lambda T^{*}M}\ e(x,\lambda )dx,
\lambda\in {\Bbb R}.
\end{equation}

By \cite{tang}, for any $\lambda  \in {\Bbb R}$,
the function
$\hbox{Tr}_{\Lambda T^{*}M}\ e(x,\lambda )$ is a
bounded measurable function on $M$, therefore,
the spectrum distribution function $N_{\cal F}(\lambda)$
is well-defined and takes finite values.

\begin{theorem}
\label{fun}
For any function $f\in C^{\infty}(\Bbb R)$, we have the following asymptotic
formula:
\begin{equation}
\label{(0.3)}
\hbox{tr}\ f(L_{h}) = h^{-q}
\frac{(4\pi)^{-q/2}}{\Gamma(q/2)}
\int_{-\infty}^{\infty}\int_{-\infty}^{\infty} \sigma^{q/2-1} f(\tau+\sigma)\
d\sigma\ dN_{\cal F}(\tau ) + O(h^{1-q}), h\rightarrow 0.
\end{equation}
\end{theorem}
\noindent{\bf Proof.} Let $E_{g_H}(\tau)$ and $E_{\Delta}(\sigma)$
denote the spectral projections of the operators
$g_H$ and $\Delta_{{\cal F}_H}$ in $L^2(H^{*},\pi^{*}\Lambda T^{*}M)$. Then,
since these operators commute, we
have
\begin{eqnarray*}
f(\sigma_h(\Delta_h)) & = & f(\Delta_F+g_H)\\
& = & \int^{+\infty}_{-\infty} \int^{+\infty}_{-\infty}
f(\tau+\sigma)\ dE_{g_H}(\tau)\ dE_{\Delta}(\sigma)\\
\end{eqnarray*}
is a tangential operator on $H^*$ with respect to the
foliation ${\cal F}_H$, which tangential kernel has the form
\begin{equation}
k_{f(\sigma_h(\Delta_h))}(\gamma,\nu,\eta)
= \int^{+\infty}_{-\infty} \int^{+\infty}_{-\infty}
f(\tau+\sigma)\ dE_{g_H}(\tau)(\nu)\ dE_{\Delta}(\gamma,\sigma).
\end{equation}
So we obtain
\begin{eqnarray}
tr_{{\cal F}_H} f(\sigma_h(\Delta_h)) & = &
\int_M \int_{H^{*}_x}\hbox{Tr}_{\Lambda T^{*}M}\
k_{f(\sigma_h(\Delta_h))}(x,\nu)dxd\nu\nonumber\\
&=&\int_M \int^{+\infty}_{-\infty} \int^{+\infty}_{-\infty}
f(\tau+\sigma)\ (\int_{H^{*}_x}dE_{g_H}(\tau)(\nu)\ d\nu)
\nonumber\\
& & d_{\sigma}(\hbox{Tr}_{\Lambda T^{*}M}\ E_{\Delta}(x,\sigma))
\ d\tau\ dx,
\end{eqnarray}
from where, taking into account that
\begin{equation}
E_{g_H}(\tau)(\nu)=\chi_{\{g_H(\nu)\leq \tau\}}I_{\pi^{*}\Lambda T^{*}M}
\end{equation}
and
\begin{equation}
\int_{H^{*}_x}E_{g_H}(\tau)(\nu)\ d\nu =
volume\{\nu\in H^{*} : g_H(\nu)\leq \tau\} =
\omega_q\tau^{q/2},
\end{equation}
where
\begin{equation}
\omega_q=\frac{\pi^{q/2}}{\Gamma((q/2)+1)}
\end{equation}
is the volume of the unit ball in ${\Bbb R}^q$,
we immediately obtain the desired formula.
\medskip
\par
In a particular case when $f$ is a characteristic function
of the semiaxis $(-\infty,\lambda)$, Theoremla~\ref{fun}
gives the asymptotic formula for the spectrum
distribution function $N_h(\lambda)$.

\begin{theorem}
\label{mainN}
Under current hypothesis, we have
\begin{equation}
N_h(\lambda) = h^{-q} \frac{(4\pi)^{-q/2}}{\Gamma((q/2)+1)}
\int_{-\infty}^{\lambda}
(\lambda - \tau )^{q/2}\ dN_{\cal F}(\tau ) + o(h^{-q}), h\rightarrow 0
\end{equation}
for any $\lambda \in {\Bbb R}$.
\end{theorem}

Theorem~\ref{intr} is, exactly, Theorem~\ref{mainN}
formulated in terms of the operator $\Delta_h$.

\section{Limits of eigenvalues}
\label{limits}
Here we discuss the asymptotical behaviour of individual
eigenvalues of the operator $\Delta_h$ when $h$ tends
to zero.

As usual, we will, equivalently, consider the operator
$L_h$ instead $\Delta_h$. Moreover, we will consider
eigenvalues of this operator on differential $k$-forms.
Therefore, we will write $L^k_h$ for the restriction
of the operator $L_h$ on $C^{\infty}(M,\Lambda^kT^{*}M)$
$k=1,\ldots,n$, omitting $k$ where it is not essential.

For any $h>0$, $L_h$ is an analytic family of type (B)
of self-adjoint operators in sence of \cite{Kato}.
Therefore, for $h>0$, the eigenvalues of $L_h$ depends
analytically on $h$. Thus there are (countably many)
analytic functions $\lambda_i(h)$ such that
\begin{equation}
\hbox{spec}\ L_h =\{\lambda_i(h):i=1,2,\ldots\},
 h>0.
\end{equation}
Moreover, by \cite{Kato}, the functions $\lambda_i(h)$
satisfy the following equality
\begin{equation}
\label{deriv}
\lambda_i'(h)=((dL_h/dh) v_h,v_h),
\end{equation}
where $v_h$ is a normalized eigenvector associated with
the eigenvalue $\lambda_i(h)$.

\begin{proposition}
\label{lim}
Under current hypotheses, for ant $i$, there exists a
limit
\begin{equation}
\label{tlim}
\lim_{h\rightarrow 0+} \lambda_i(h)=\lambda_{\lim,i}.
\end{equation}
Moreover, if $v_h$ is a normalized eigenform associated with
the eigenvalue $\lambda_i(h)$, then we have the estimates
\begin{eqnarray}
\label{eigen}
\|v_h\|_{0,1}<C_1,\ h\|v_h\|_{1,0}<C_2,
\end{eqnarray}
with constants $C_1$ and $C_2$ independent of $h\in(0,1]$.
\end{proposition}
\noindent{\bf Proof.}\
As above, let $v_h$ be a normalized eigenform associated with
the eigenvalue $\lambda_i(h)$:
\begin{equation}
L_hv_h=\lambda_i(h)v_h,\ \|v_h\|=1.
\end{equation}
By (\ref{deriv}), we have
\begin{equation}
\lambda_i'(h) = ((2h\Delta_H + 4h^3\Delta_{-1,2}+
 K_1+2hK_2 +3h^2K_3)v_h, v_h),
\end{equation}
from where, using the positivity of operators
$\Delta_H$ and $\Delta_{-1,2}$ in $L^2(M,\Lambda T^{*}M)$,
 and the estimates
(\ref{K2}) and (\ref{K1}) (with $h=1$), we obtain
\begin{equation}
\label{lam}
\lambda_i'(h)\geq -C_1\|v_h\|^2_{0,1}-C_2h^2\|v_h\|^2_{1,0}
-C_3.
\end{equation}
The estimate~(\ref{crude}) implies
\begin{equation}
\label{est}
C_1\|v_h\|^2_{0,1}+C_2h^2\|v_h\|^2_{1,0}\leq C_3\lambda_i(h)+
C_4, h\in (0,1].
\end{equation}
By~(\ref{lam}) and (\ref{est}), we conclude that
\begin{equation}
\lambda_i'(h)\geq -C_5\lambda_i(h)-C_6.
\end{equation}
This estimate can be rewritten in the following way:
\begin{equation}
\frac{d}{dh}((\lambda_i(h)+ \frac{C_6}{C_5})e^{C_5h})\geq 0,
\end{equation}
that means that the function
 $(\lambda_i(h)+ \frac{C_6}{C_5})e^{C_5h}$ is increasing
in $h$ for $h$ small enough. By the positivity of
the operator $L_h$ in $L^2(M,\Lambda T^{*}M)$, every
eigenvalue $\lambda_i(h)$ is positive, so the function
 $(\lambda_i(h)+ \frac{C_6}{C_5})e^{C_5h}$
semibounded from below near zero. Therefore, this function
has a limit when $h$ tends to zero, that, clearly, implies
the existence of the limit for the function
$\lambda_i$.

The second assertion of this Proposition is an immediate
consequence of the first one and the estimate~(\ref{est}).
\medskip
\par
Proposition~\ref{lim} allows us to introduce the limitting spectrum of the
operator $\Delta^k_h$ as a set of all limitting values
$\lambda^k_{\lim,i}$, given by (\ref{tlim}):
\begin{equation}
\sigma_{\lim}( \Delta^k_h)=\{\lambda^k_{\lim,i}:
i=0,1,\ldots\}.
\end{equation}

By an analogy with the case of semiclassical asymptotics
for Schrodinger operator, we may assume that
the structure of the
limitting spectrum $\sigma_{\lim}(\Delta^k_h)$ is defined
in a big extent by a limitting
value of the bottoms of spectrum of the operator $\Delta^k_h$. So let
\begin{equation}
\lambda^k_0(h)=\min_{u\in C^{\infty}(M,\Lambda^k
T^{*}M)}\frac{(\Delta^k_hu,u)}{\|u\|^2},
\end{equation}
and
\begin{equation}
\lambda^k_{\lim,0}=\lim_{h\rightarrow 0}\lambda^k_0(h).
\end{equation}
There are two other quantities: the bottom $\lambda^k_{F,0}$
of the spectrum of the operator $\Delta^k_F$ in
$L^2(M,\Lambda^k T^{*}M)$:
\begin{equation}
\lambda^k_{F,0}=
\min_{u\in C^{\infty}(M,\Lambda^k T^{*}M)}
\frac{(\Delta^k_Fu,u)}{\|u\|^2},
\end{equation}
and the bottom $\lambda^k_{{\cal F},0}$
of the leafwise spectrum of the operator $\Delta^k_F$ in
$L^2(L,\Lambda^k T^{*}M)$:
\begin{equation}
\lambda^k_{{\cal F},0}=\overline{\bigcup\{\sigma(\Delta^k_L):
L\in V/{\cal F}\}},
\end{equation}
where
\begin{equation}
\lambda^k_{L,0}=
\min_{u\in C^{\infty}(L,\Lambda^k T^{*}M)}
\frac{(\Delta^k_Lu,u)}{\|u\|^2},
\end{equation}
the operator $\Delta^k_L$ is the restriction of the operator
$\Delta^k_F$ on the leaf $L$.

\begin{proposition}
Under current hypotheses, we have the following relations:
\begin{equation}
\label{rel}
\lambda^k_{F,0}\leq \lambda^k_{\lim,0}\leq  \lambda^k_{{\cal F},0},\
k=1,\ldots,n.
\end{equation}
\end{proposition}
\noindent{\bf Proof.}\ Let $v_h$ be a normalized eigenform associated with
the bottom eigenvalue $\lambda^k_0(h)$:
\begin{equation}
L^k_hv_h=\lambda^k_0(h)v_h,\ \|v_h\|=1.
\end{equation}
By the definition of $\lambda^k_{F,0}$, we have the estimate
\begin{equation}
(\Delta^k v_h,v_h)\geq \lambda^k_{F,0}.
\end{equation}
By (\ref{Lh3}), we obtain
\begin{equation}
\label{Lh4}
\lambda^k_0(h)\geq (1-h^2)\lambda^k_{F,0}+C_1h^2\|v_h\|^2_{1,0}
h(K_1v_h,v_h)- C_2h^2,
\end{equation}
where $C_1$ and $C_2$ are positive constants.
By (\ref{eigen}), we have
\begin{equation}
\lim_{h\rightarrow 0}h(K_1v_h,v_h)=0.
\end{equation}
Taking this into account, by (\ref{Lh4}), we immediately complete
the proof of the inequality
$\lambda^k_{F,0}\leq \lambda^k_{\lim,0}$.

Theorem~\ref{intr} implies that
$N^k_h(\lambda)>0$ for any $\lambda>\lambda^k_{{\cal F},0}$
and $h$ small enough, from where the desired inequality $\lambda^k_{\lim,0}\leq
 \lambda^k_{{\cal F},0}$
follows immediately.
\medskip
\par
We conclude this Section with some remarks and examples,
concerning quantities $\lambda^k_{F,0}$, $\lambda^k_{\lim,0}$
and $\lambda^k_{{\cal F},0}$.
\medskip
\par
\noindent{\bf Remarks.}\ (1)
When the foliation ${\cal F}$ is a fibration or, more
general, ${\cal F}$ is amenable in some sense (see also
Section~\ref{disc}),
relations~(\ref{rel}) turns out to be identities \cite{tang}.
\medskip
\par
\noindent (2) We don't know if the equality
$\lambda^k_{F,0}=\lambda^k_{\lim,0}$ is always true. It is, clearly, so for
$k=0$:
\begin{equation}
\lambda^0_{F,0}=\lambda^0_{\lim,0}=0.
\end{equation}
Another remark is as follows. If the Betti number $b_k(M)$ is not
zero, then $\lambda^k_0(h)=0$ for all $h$, that also implies
$$
\lambda^k_{F,0}=\lambda^k_{\lim,0}=0.
$$
\noindent (3) \ Here we give an example of the
foliation such that the bottom $\lambda^0_{F,0}=0$
of the operator $\Delta^0_F$ in $L^2(M)$ is a point of discrete spectrum.
\medskip
\par
\noindent {\bf Example.}\ Let $\Gamma$ be a discrete, finitely generated group
such that

(a)\ $\Gamma$ has property $(T)$ of Kazhdan;

(b)\ $\Gamma$ is be embedded in a compact Lie group $G$ as a dense subgroup.

\noindent For definitions and examples of such groups, see, for
instance,\cite{JV,Lub}.

Let us take a compact manifold $X$ such that
$\pi_1(X)=\Gamma$. Let $\tilde{X}$ be the universal
covering of $X$ equipped with a left action of $\Gamma$
by deck transformations. We will assume that $\Gamma$
acts on $G$ by left translations.
Let us consider the suspension foliation
${\cal F}$ on a compact manifold $M=\tilde{X}\times_{\Gamma}G$ (see, for
instance, \cite{Cam}). A choice of a left invariant
metric on $G$ provides a bundle-like metric on $M$,
so ${\cal F}$ is a Riemannian foliation.
We may assume that leafwise
metric is chosen in such a way that any leaf of the
 foliation ${\cal F}$ is isometric to $\tilde{X}$.

There is defined a natural action of $\Gamma$ on $M$
and the operator $\Delta^0_F$ is invariant under
this action. Let $E(0,\lambda),\lambda>0,$ denote the spectral
projection of the operator $\Delta^0_F$ in $L^2(M)$,
corresponding to the interval $(0,\Lambda)$, and
$E(0,\lambda)L^2(M)$ be the corresponding $\Gamma$-invariant
spectral subspace.
\medskip
\par
\noindent {\it Claim.}\ In this example, the bottom $\lambda^0_{F,0}=0$
of leafwise Laplacian in $L^2(M)$ is a nondegenerate
point of discrete spectrum of the operator
$\Delta^0_F$, that is, an isolated eigenvalue of the
multitplicity 1.
\medskip
\par
{}From the contrary, let us assume that zero lies in the
essential spectrum of the operator $\Delta^0_F$ in
$L^2(M)$. Then, for any $\varepsilon>0$ and $\lambda>0$, there is
a function $u_{\varepsilon}\in C^{\infty}(M)$ such that
$u_{\varepsilon}$ belongs
to the space
$E(0,\lambda)L^2(M)$,  $\|u_{\varepsilon}\|=1$ and
\begin{equation}
\label{df}
(\Delta_Fu_{\varepsilon}, u_{\varepsilon})
=\|\nabla_Fu_{\varepsilon}\|\leq \varepsilon,
\end{equation}
where $\nabla_F$ denotes the leafwise gradient. From (\ref{df}),
we can easily derive that the representation of
the group $\Gamma$ in $E(0,\lambda)L^2(M)$ has almost
invariant vector, that, by the property $(T)$, implies
the existence of an invariant vector
$v_0\in E(0,\lambda)L^2(M)$.

Since $\Gamma$ is dense in $G$, $\Gamma$-invariance
of $v_0$ implies its $G$-invariance, that, in turn,
implies that $v_0$ is a lift of some non-zero element
$v\in C^{\infty}(X)$
via the natural projection $M\rightarrow X$. It can be
easily checked that $v$ belongs to the corresponding
spectral space $E(0,\lambda)L^2(X)$ of the Laplace
operator $\Delta_X$ in $L^2(X)$. From other hand,
the operator $\Delta_X$
has a discrete spectrum, so zero is an isolated point
in the spectrum of $\Delta_X$, and $E(0,\lambda)L^2(X)$
is a trivial space if $\lambda>0$ is small enough.
So we get a contradiction, which imply that zero lies in the
discrete spectrum of the operator $\Delta^0_F$ in
$L^2(M)$.
\medskip
\par
\noindent(4)\ In the case of a fibration, we also have that zero is
an isolated point in the
spectrum of the operator $\Delta^0_F$ in
$L^2(M)$, but, in this case, it is an eigenvalue
of infinity multitplicity, so that it lies in the
essential spectrum of the operator $\Delta^0_F$ in
$L^2(M)$.
\medskip
\par
\noindent (5) \ Unlike the scalar case, it is not always the case that all of
the semiaxis $[\lambda_{\lim,0},+\infty)$ is contained
in $\sigma_{lim}(\Delta_h)$.
Indeed, let, as in
Example of (3), $\lambda^0_{F,0}=0$
is a nondegenerate point of discrete spectrum of the operator $\Delta^0_F$.
Then, by means of the perturbation theory of the discrete spectrum (see, for
instance,
\cite{Kato}), we can state that, for $h>0$ small enough,  $\lambda^0(h)=0$ is
the only eigenvalue of the operator $\Delta^0_h$
near zero. So we conclude that there exists a
$\lambda_1>0$ such that, for any $h>0$ small enough,
\begin{equation}
\sigma_{lim}(\Delta_h)\bigcap [\lambda_1,+\infty)=0.
\end{equation}

\section{Some remarks on the main asymptotical
formula}
\label{disc}
In this Section, we discuss
some aspects of the main
asymptotical formula~(\ref{eig1}). We are, especially,
interested in a discussion of the formula~(\ref{eig2}).
We will make use of the notation of previous Sections.

So recall that the whole picture which we observe in  the foliation
case is the following.
In a general case, for any $k=0,1,\ldots,n$, we have only that
\begin{equation}
\lambda^k_{F,0}\leq \lambda^k_{\lim,0}\leq  \lambda^k_{{\cal F},0},
\end{equation}
and these relations turns into identities, if
the foliation ${\cal F}$ is a fibration or, more
general, is amenable in some sense (see Section~\ref{limits}
and \cite{tang} for discussion).

By (\ref{eig1}), the function $N^k_h(\lambda)$
behaves as usual when $\lambda$ is greater than the bottom
of the leafwise spectrum of $\Delta^k_F$:
\begin{equation}
N^k_h(\lambda)\sim Ch^{-q}, \lambda \geq \lambda^k _{{\cal F},0},
\end{equation}
but, if $\lambda^k_{F,0}<\lambda^k_{{\cal F},0}$,
there might be limitting
values for eigenvalues $\lambda^k_i(h)$ of the operator
$\Delta^k_h$, lying in the interval
$(\lambda^k_{F,0}, \lambda^k _{{\cal F},0})$.
So the function $N^k_h(\lambda)$ is nontrivial on the
interval
$(\lambda^k_{\lim,0},
\lambda^k _{{\cal F},0})$,
but the fact mentioned above that the right-hand side of (\ref{eig1})
 depends only
on leafwise spectral data of the operator $\Delta^k_F$
 implies the formula
\begin{equation}
\lim_{h\rightarrow 0+}h^{q}N^k_{h}(\lambda ) = 0,\ \lambda  < \lambda^k _{{\cal
F},0}.
\label{(0.8)}
\end{equation}
It means that the set of eigenvalues of
$\Delta^k_{h}$ in the interval $(\lambda^k_{\lim,0},
\lambda^k _{{\cal F},0})$ is "thin" in the whole set of eigenvalues of
$\Delta_{h}$.
By analogy with \cite{Su}, (\ref{(0.8)}) in the case
$k=0$ may be called as a weak foliated version of "Riemann hypothesis".

This is quite different from what we have in the
case of Schrodinger operator or in the fibration case.
For instance, if $H_h$ is
the Schrodinger operator on a compact manifold $M$
(we may consider $M$, being equipped with a trivial
foliation ${\cal F}$ which leaves are points):
\begin{equation}
H_h=-h^2\Delta +V(x), x\in M.
\end{equation}
we have
\begin{equation}
\lambda_{F,0}= \lambda_{\lim,0}=\lambda_{{\cal F},0}
=\inf V_{-},
\end{equation}
where
\begin{equation}
V_{-}(x)=\min(V(x),0), x \in M,
\end{equation}
and the following asymptotical formula for spectrum distribution
function $N_h(\lambda)$ in semiclassical limit:
\begin{equation}
\label{semi}
N_h(\lambda)=(2\pi)^{-n}h^{-n}\int_{\{(x,\xi): \xi^2+V(x)
\leq\lambda\}}\ dxd\xi+o(h^{-n}), h\rightarrow 0+.
\end{equation}
So, if $h\rightarrow 0$. the picture is as follows:
\begin{equation}
N_h(\lambda)\sim Ch^{-n}, \lambda>\inf V_{-},
\end{equation}
where $n=\dim M$ and
\begin{equation}
N_h(\lambda)=0, \lambda\leq\inf V_{-}.
\end{equation}

It is worthwhile to note facts in spectral theory of coverings, which are very
similar to
ones in spectral theory of foliations mentioned above.
Let us consider the case of Laplace-Beltrami
operator on functions.

Let $\tilde{M}\rightarrow M$ be a normal covering with
a covering group $\Gamma$. Recall that a tower of coverings is a set
$\{M_i\}_{i=1}^{\infty}$ of finite-fold subcoverings
of this covering with the corresponding covering groups
$\Gamma_i$ such that:

(1) for each $i$, $\Gamma_i$ is a normal subgroup of
finite index in $\Gamma$;

(2) for each $i$, $\Gamma_{i+1}$ is contained in $\Gamma_i$;

(3) $\bigcap_i\Gamma_i = \{e\}$.

Let $\sigma(\Delta_{M_i})$ be a set of eigenvalues of the
Laplacian on $M_i$, and $N_{M_i}(\lambda)$ be its
distribution function.
For
any $i$, we have an embedding
\begin{equation}
\label{emb}
\sigma(\Delta_{M_i})\subset\sigma(\Delta_{M_{i+1}}),
\end{equation}
and when $i$ tends to the infinity
the spectrum $\sigma(\Delta_{M_i})$ of a finite covering
approaches to a limit
\begin{equation}
\sigma_{lim}(\Delta)=\bigcup_i\sigma(\Delta_{M_i}).
\end{equation}
Then, the bottom $\lambda_{\lim,0}$ of limitting spectra $\sigma _{\lim
}(\Delta)$ and the bottom $\lambda_{M,0}$ of the spectrum
$\sigma (\Delta_M)$  of the manifold $M$ are, clearly, equal to 0.
By \cite{Br}, the bottom $\lambda_{\tilde{M},0}$ of the spectrum
$\sigma (\Delta_{\tilde{M}})$  of the covering manifold is
equal to $\lambda_{M,0}$:
\begin{equation}
\lambda_{\tilde{M},0}=\lambda_{M,0},
\end{equation}
if and only if the  group $\Gamma$ is amenable.

Moreover, by \cite{Do}, for any function $f\in C^{\infty}_c({\Bbb R})$,
we have
\begin{equation}
\lim_{i\rightarrow \infty}(vol\ M_i)^{-1}tr\ f(\Delta_{M_i})
=tr_{\Gamma}\ f(\Delta_{M}),
\end{equation}
where $tr_{\Gamma}$ is von Neumann trace on the
the algebra of $\Gamma$-invariant operators on $\tilde{M}$
\cite{At}.
In particular, if $N_{i}(\lambda)$ is the eigenvalue distribution
function of the Laplace-Beltrami operator $\Delta_{M_i}$,
then
\begin{eqnarray}
\lim_{i\rightarrow \infty}(vol\ M_i)^{-1}N_{i}(\lambda)&=&
N_{\Gamma}(\lambda),\lambda\in{\Bbb R},\\
\lim_{i\rightarrow \infty}(vol\ M_i)^{-1}N_{i}(\lambda)&=&
0,\lambda<\lambda_{\tilde{M},0},
\end{eqnarray}
where $N_{\Gamma}(\lambda)$ is spectrum distribution
function of the operator $\Delta_{\tilde{M}}$ constructed
by means of the $\Gamma$-trace $tr_{\Gamma}$,
$\lambda_{\tilde{M},0}=\inf \sigma (\Delta_{\tilde{M}})$

A little bit more general possibility to arrange finite-dimensional
approximation of the spectrum of a covering, making use of sequences of
finite-dimensional representations
of a covering group $\Gamma$, converging to the left
regular representations of $\Gamma$, is considered
in \cite{Su}.
Analogues of (\ref{eig1}) and (\ref{(0.8)}) can be also
found in \cite{Su}.

We may point out two common features
of spectral theory
for Laplacian on a covering and spectral theory for leafwise Laplacian
on foliated manifold. From the tangential point of view, both of them can be
treated as type II spectral
problems in a sense of theory of operator algebras, and
asymptotical spectral problems mentioned above can be
considered as finite-dimensional (of type I) approximations to these spectral
problems. Actually, some spectral
characteristics related with such an approximation
don't depend on a choice of a bundle-like metric
on $M$, and, moreover, are invariants of quasi-isometry
of metrics (coarse invariants in a sense of
\cite{Ro}).
One of the simplest characteristics of such a kind which we have already
met is the notion of amenability.

We can introduce some quantative spectral characteristics of the tangential
Laplacian $\Delta^k_F$
related with adiabatic limits.
For any $\lambda$, let $r_k(\lambda)$ be given
as
\begin{equation}
r_k(\lambda)=-\limsup_{h\rightarrow 0}\ln N^k_h(\lambda)/\ln h.
\end{equation}
Otherwise speaking, $r_k(\lambda)$ equals
the least bound of all $r$ such that
\begin{equation}
\label{r}
N^k_h(\lambda)\sim Ch^{-r}, h\rightarrow 0.
\end{equation}
If $\lambda<\lambda^k_{\lim,0}$, we put $r_k(\lambda)=
-\infty$.

Then we can easily state the following properties of the function
$r_k(\lambda)$:
\begin{enumerate}
\item $0\leq r_k(\lambda)\leq q$ for any
$\lambda\geq \lambda^k_{\lim,0}$;
\item $r_k(\lambda)$ is not decreasing in $\lambda$;
\item  $r_k(\lambda)=q$ if $\lambda>\lambda^k _{{\cal F},0}$.
\item if the foliation ${\cal F}$ is amenable, then:
\begin{eqnarray*}
&r_k(\lambda)=q,&\lambda>\lambda^k _{{\cal F},0},\\
&r_k(\lambda)=-\infty,&\lambda\leq\lambda^k _{{\cal F},0}.
\end{eqnarray*}
\item $r_k(\lambda)=0$ iff the interval $[0,\lambda]$
lies in the discrete spectrum of the operator
$\Delta^k_F$ in $L^2(M,\Lambda^kT^{*}M)$. As we have seen
in the previous Section, such situation can happen
(a property $(T)$ case).
\end{enumerate}
Then we expect that some invariants of the function
$r_k(\lambda)$ introduced above  near $\lambda=0$ might to
be independent of the choice of
metric on $M$ (otherwise speaking, to be coarse invariants),
and, moreover, be topological or homotopic
invariants of foliated manifolds.

{}From transversal point of view, both of them are related with some sort of
"noncommutative" fibration in sense of noncommutative
differential geometry \cite{Co}. Here the relation
(\ref{(0.8)}) reflects a nontriviality of geometry of these "fibrations" in the
nonamenable case.

Now we point out two facts
in noncommutative spectral geometry of foliations,
which are closely related with (\ref{(0.8)}).
When the foliation ${\cal F}$ is Riemannian, we
can consider $M/{\cal F}$ as a noncommutative Riemannian manifold.
More precise, we can  define the corresponding
spectral triple (in a sense of \cite{Co-M}) as follows:
\begin{enumerate}
\item An involutive
algebra ${\cal A}$ is an algebra $C^{\infty}_c(G_{\cal F})$ of smooth,
compactly supported functions on the holonomy groupoid
$G_{\cal F}$ of the foliation ${\cal F}$;
\item A Hilbert space ${\cal H}$
is a space $L^2(M,\Lambda H^{*})$ of the transversal
differential forms, on which
an element $k$ of the algebra ${\cal A}$ is represented via
a smoothing tangential operator with the tangential
kernel $k$;
\item an operator $D$ is the transverse signature
operator $d_H+\delta_H$ of a bundle-like metric on $M$.
\end{enumerate}

Let $C^{*}(G_{\cal F})$ ($C^{*}_r(G_{\cal F})$) be the full (reduced)
$C^{*}$-algebra of the foliation respectively.
There is the natural
projection $\pi:C^{*}(G_{\cal F})\rightarrow C^{*}_r(G_{\cal F})$.
We say that the foliation ${\cal F}$ is amenable, if
the
projection $\pi:C^{*}(G_{\cal F})\rightarrow C^{*}_r(G_{\cal F})$ is an
isomorphism.

The first fact is that, in a case of the foliation ${\cal F}$ is nonamenable,
this noncommutative Riemannian manifold has pieces
of various dimension with the top dimension, being,
certainly, equal to $q$ in the following sense.

Let us consider subsets of $V/{\cal F}$ as involutive
ideals in $C^{*}(G_{\cal F})$.
We can speak about the top
spectral dimension of the pieces of our space which are
contained ${\cal I}$ in the following way (see \cite{Co-M}
for details).
We say that this bound is less than $k$, if for any $a\in {\cal I}$ the
distributional zeta function
\begin{equation}
\zeta_a(z)=\hbox{tr}\ a|D|^{-z}
\end{equation}
 extends
holomorphically to the halfplane $\{z\in {\Bbb C}:
Re\ z >k\}$.
By the Tauberian theorem,
the top dimension of the subset in the space can be also
detected by means of asymptotics of the distributional spectrum distribution
function
\begin{equation}
N_a(\lambda)=tr(aE_{\lambda}(|D|)), a\in {\cal I},
\lambda\in {\Bbb R},
\end{equation}
where $E_{\lambda}(|D|)$ is the spectral projection
of the operator $|D|$, corresponding to the semiaxis
$(-\infty,\lambda)$, or the theta-function
\begin{equation}
\theta_a(t)=tr(ae^{-tD^2}), a\in {\cal I}, t>0.
\end{equation}
For instance, the top
spectral dimensions of the pieces of our space which are
contained ${\cal I}$ is less than $k$, if for any
$a\in {\cal I}$ the
distributional theta function $\theta_a(t)$ satisfies the
estimate
\begin{equation}
\label{top}
\theta_a(t)\leq Ct^{-k/2}, 0<t\leq 1.
\end{equation}

Then we have (compare with Proposition 4.4 in \cite{asymp}):
\medskip
\par
An involutive ideal ${\cal I}$
in $C^{*}(G_{\cal F})$ has the top dimension $q$ iff ${\cal I}\bigcap
\pi(C^{*}_r(G_{\cal F}))\not=\emptyset$.
In particular, if
$\pi({\cal I})=0$, then the top spectral dimension
of ${\cal I}$ is less than  $q$.
\medskip
\par
The other fact is related with
the support of the "noncommutative" integral,
given by the Dixmier trace $Tr_{\omega}$. Namely, it can
be shown that in the case under consideration
the Dixmier trace $Tr_{\omega}(k)$, corresponding to
the spectral triple introduced above exists and doesn't
depend on a choice of $\omega$ for any $k\in C^{*}(G_{\cal F})$.
Then we have
\begin{equation}
Tr_{\omega}(k)=0
\end{equation}
for any $k\in C^{*}(G_{\cal F}),\pi(k)=0$
{}.
To relate these facts with the spectral theory
of the tangential Laplace-Beltrami operator $\Delta_F$,
we have to
note that, by \cite{tang}:
\medskip
\par
(1) the operator $f(\Delta_F)$
belongs to the $C^{*}$- algebra
$C^{*}(G_{\cal F})$ for any
$f\in C^{\infty}_c({\Bbb R})$, and,

(2) by spectral theory, $\pi(f(\Delta_F))=0$
for any $f\in C^{\infty}_c({\Bbb R})$ such that
$supp \ f \subset (\lambda_{\lim,0},
\lambda _{{\cal F},0})$.
\medskip
\par
It seems also to be true that the function $r_k(\lambda)$
introduced above takes values in the spectrum dimension
$Sd$ of the noncommutative spectrum space in question
(see \cite{Co-M}).

\end{document}